\documentclass[11pt]{amsart}
\usepackage[english]{babel}
\usepackage[dvips]{graphicx}
\renewcommand{\div}{\mbox{\rm div}\ }
\newcommand{\mes}{\mbox{\rm mes}\ }
\newcommand{\supp}{\mbox{\rm supp}\ }
\theoremstyle{plain}
\newtheorem{thrm}{Theorem}[section]
\newtheorem{lmm}[thrm]{Lemma}

\newtheorem{prpstn}[thrm]{Proposition}

\theoremstyle{definition}
\newtheorem{dfntn}[thrm]{Definition}

\newtheorem{rmrk}[thrm]{Remark}

  \def\xR{\mathbb{R}}

\def\xCzero{{\rm C}^{0}}
\def\xCone{{\rm C}^{1}} 
 
\def\xCinfty{{\rm C}^{\infty}} 

\def\xdif{\,{\rm d}}
\usepackage[dvips]{hyperref}

\numberwithin{equation}{section}

\begin{document}
\title{On shape optimization and the Pompeiu problem}
\author{Ar\=unas Grigelionis}
\address{Inst. of Mathematics and Informatics,
Akademijos 4, 2600 Vilnius, Lithuania }
\date{}
\begin{abstract}
The Pompeiu problem is considered as shape optimization problem.
We show stability of the ball which is the minimum point of related domain functional.
The proof is based on shape derivative method. Stability of the ball for
general domain functionals invariant under the rigid motions  is discussed.
\end{abstract}
\keywords{shape optimization, positive definite functions, the Pompeiu problem\\E-mail adress: agrig@ktl.mii.lt}
\maketitle
\section*{Introduction and statements}

Energy--type functionals $\it{F}\left[ \omega \right] =\iint
f\left( \left| x-y\right| \right) \xdif {\omega (x)}\xdif {\omega (y)}$, where $\omega$ is a measure
of compact support, appears in statistical mechanics of systems, topological classification problems
of knots,  isoperimetric problems, harmonic analysis, discrete energy problems and other
areas of pure and applied mathematics (see, for example, \cite{Burch, Fri, Sjolin,  BFHW,  RSZ} ).

Motivated by problems of shape optimal design \cite{Cea, Zol, Ker, Mur, Grig1}
we consider the model domain functional

\begin{equation}\label{funk}
\it{F}\left[ \Omega \right] =\int\limits_\Omega \int\limits_\Omega
f\left( \left| x-y\right| \right) \xdif x\xdif y,
\end{equation}
defined on the set of bounded domains $\Omega \subset 
\xR^n,\quad n\geq 2,$ with smooth boundaries. Here $f:\left( 0,\infty \right) \rightarrow 
\xR$ is known function.

Let $V\in \xCone\left( \xR^n,\xR^n\right) .$ Equation
\begin{equation}\label{vect}
\dot x=V\left( x\right)
\end{equation}
will produce a flow $T_t^V:\xR \times \xR^n\rightarrow \xR^n$ which moves a bounded
smooth domain $\Omega=\Omega _0$ to its new position $\Omega
_t^V=T_t^V\left( \Omega _0\right)$ (see Fig.~\ref{fig1} on p.~\pageref{fig1}). 
For all $t$ close enough to zero the domain $\Omega _t^V$ has been shown to be bounded 
and smooth \cite{Zol}.

Shape derivative in the direction of vector field $V$ is defined by formula
\begin{equation}\label{shder}
\it{\dot F}\left[ \Omega _0;V\right] :=\frac d{dt} |_{t=0}\it{F}\left[ \Omega _t^V\right] .   
\end{equation}

\begin{dfntn}\label{critic} 
We say nonempty domain $\Omega$ is \textit{critical} if 
$$
\it{\dot F}\left[ \Omega ;V\right]=0 
$$
for any $V\in \xCone\left( \xR^n, \xR^n\right) $. 
\end{dfntn}
If $\Omega$ is critical, then (see Section \ref{shpdrv}) equality
\begin{equation} \label {eq0}
\int\limits_\Omega f\left( \left| x-y\right| \right) \xdif y=const
\end{equation}
holds for all $x\in \partial \Omega $ with $const=0$.

As a special case equation (\ref{eq0}) includes the inverse potential problem
\begin{equation} \label {eq01}
\int\limits_\Omega \left| x-y\right| ^{2-n} \xdif y=const.
\end{equation}
This last equation defines the body which has constant gravity potential on its own shape. 
The observation based on the "moving planes" method was done in \cite{Grig2}:
\begin{prpstn}\label{monot}
Let $\xCone$ ~function $f $ is strictly monotone. If nonempty $\xCone$ domain $\Omega$  
solves equation (\ref{eq0}), then it is a ball.
\end{prpstn}
Thus, if $\xCone$ ~function $f $ is strictly monotone, then any critical point of the domain functional (\ref{funk}) 
is the ball of fixed radius. If function $f$ is not strictly monotone, then critical domains of different shapes 
are possible. 

Usually {\it stable} minimum point of general domain functional $\it{F}$ is 
defined as the domain $\Omega_0$ at which inequality $\it{F}\left[ \Omega_0 \right]<\it{F}\left[ \Omega\right]$
holds for all $\Omega$ close enough to $\Omega_0$ (see, for example, \cite{DP}). Even despite of the
difficulty in proper understanding of the closeness, domain functional (\ref{funk}) has no stable minimum
(and maximum) points in this sense because it is invariant under the rigid motions. 

To reach the uniformity with definition of critical point we shall introduce the notion of {\it weakly stable} minimum point
for the domain functional $\it{F}$.
\begin{dfntn}\label{defstb}
We say the domain $\Omega=\Omega_0$ is stable minimum point of the domain functional (\ref{funk}) in 
the direction of vector field $V\in \xCone\left( \xR^n,\xR^n\right)$, if 

either for all $t$, $t\neq 0$ and close enough to zero, inequality
$\it{F}\left[ \Omega_0 \right]<\it{F}\left[ \Omega _t^V\right]$
holds,

either for all $t$ close enough to zero,  $\Omega_0$ as the rigid body coincides with $\Omega _t^V$.
\bigskip

The domain $\Omega=\Omega_0$ is said to be {\it weakly stable} minimum point of the domain functional (\ref{funk}) in 
$\xCone\left( \xR^n,\xR^n\right)$, if 
this domain is stable minimum point of $\it{F}$ in 
the direction of vector field $V\in \xCone\left( \xR^n,\xR^n\right)$
for any vector field $V\in \xCone\left( \xR^n,\xR^n\right)$. 
\end{dfntn}
\begin{rmrk}\label{rmrk1}
{\bf a)} Inequality
\begin{equation} 
\it{\ddot F}\left[ \Omega_0;V,V\right] > 0,
\end{equation}
here $\it{\ddot F}$ denotes the second shape derivative
of the domain functional $\it{F}$, is sufficient for the critical point $\Omega=\Omega_0$ to be stable minimum
in the direction of vector field $V$.

{\bf b)} If $V=const$, then any $\Omega_0$ as the rigid body coincides with $\Omega_t^V$ for all $t\in \xR$.

{\bf c)} If normal part of $V$ on the boundary $\partial \Omega_0$ is zero, then $\Omega_0 = \Omega_t^V$
for all $t$ close to zero.
\end{rmrk}

Statements {\bf b)} and {\bf c)} are easy to understand from geometrical point of view. They follow from the standard
properties of ordinary differential equations.
\bigskip

If $\xCone$ ~function $f $ is strictly monotone, then by Proposition \ref{monot} the domain functional $\it F$
has at most one critical point to be considered as a rigid body which, of course, is weakly stable in $\xCone\left( \xR^n,\xR^n\right)$.

In this paper we consider domain functional (\ref{funk}) when  function  
$f\left( \left| \cdot \right| \right)$ is positive definite, that is
\begin{equation} \label {ineq1}
\iint f\left( \left| x-y\right| \right) \mu \left( x\right) \Bar{\mu} \left( y\right) \xdif x\xdif y\geq 0
\end{equation}
for all complex $\mu \in \xCinfty_0\left( \xR^n\right)$ and thus inequality
\begin{equation}\label{ineq2}
\it{F}\left[ \Omega \right] \geq 0
\end{equation}
holds for all  $\Omega $ bounded.

\begin{thrm}\label{stbl}
Let $\it{F}$ is the domain functional (\ref{funk}) generated by positive definite function $f\left( \left| \cdot \right| \right)$ and
\begin{equation}\label{eq1}
\it{F}\left[ \Omega\right]=0
\end{equation} 
for some $\Omega$ the ball.

Then $\Omega$ is weakly stable minimum point of $\it{F}$ in $\xCone\left( \xR^n,\xR^n\right)$.
\end{thrm}
Investigation of the domain functional (\ref{funk}) with the positive definite function $f\left( \left| \cdot \right| \right)$ is stimulated because of  its relation to the Pompeiu problem \cite{Zalc}. 
We will prove
\begin{prpstn}\label{pompeiu} Let $\Omega$ be bounded nonempty domain. 

{\bf a)} If $\Omega$ satisfies 
(\ref{eq1})  for some positive definite function
$f\left( \left| \cdot \right| \right)$, then there is $\lambda >0$:
\begin{equation} \label{eq2}
\hat \chi _{\Omega}\left( \xi \right) =0 
\quad\text{for all}\quad \xi \in \xR^n,\quad\left| \xi \right| =\lambda,
\end{equation}
here $\hat \chi _{\Omega}$ denotes the Fourier transform of indicator function $\chi _{\Omega}$.
\bigskip

{\bf b)} On the other hand, if condition (\ref{eq2}) holds, then $\it{F}\left[ \Omega \right] =0$ for 
some positive definite function $f\left( \left| \cdot \right| \right)$.
\end{prpstn}
Domain $\Omega \subset \xR^n$ which satisfies (\ref{eq2}) is known as domain without the {\it Pompeiu property}
\cite{BST}. If  its boundary $\partial \Omega$ is smooth and connected, then overdetermined problem
\begin{subequations}\label{schif}
\begin{eqnarray}
\Delta u+\lambda ^2u=0 &\mbox{in}&\Omega \\
u=1 &\mbox{on}&\partial \Omega\\
\nabla u=0 &\mbox{on}&\partial \Omega
\end{eqnarray}
\end{subequations}
has nontrivial solution \cite{Ber1, Will}. It is expected  (Schiffer conjecture \cite{Yau}) such the domain to be necessarily a ball.

In \cite{AgS, Kob} the Schiffer conjecture is affirmed for "small" perturbations of the ball.
Kobayashi \cite{Kob} use precise Fourier transform estimates of domain indicator function.
Agranovsky and Semenov \cite{AgS}, instead, exploit the overdetermined problem (\ref{schif}).
Their approach leads to the generalization of results in Riemannian spaces. 
More on the Pompeiu and Schiffer problems see the expository paper \cite {Ber2}.

The proof of the Theorem \ref{stbl} may be derived from \cite {AgS, Kob}. Instead of, we consider the Pompeiu problem as
shape optimization problem and shall give the direct proof of the Theorem \ref{stbl} based on the shape derivative method.

The paper is organized as follows. In Section \ref{shpdrv}  shape derivatives of the domain functional (\ref{funk})
are calculated. 
The domain equation (\ref{eq0}) which describes critical domains is derived and an expression of the 
second shape derivative of the domain functional (\ref{funk}) in a case of positive definite function 
$f\left( \left| \cdot \right| \right)$ is examined.

The proof of the Theorem \ref{stbl} is preceded by
\begin{lmm}\label{lem1}
Let $\it{F}$ is the domain functional (\ref{funk}) generated by positive definite function $f\left( \left| \cdot \right| \right)$ and
$\it{F}\left[ \Omega\right]=0$ for $\Omega$ the ball centered at $x=0$.
Then $\it{\ddot F}\left[ \Omega;V,V\right] = 0$ for vector field $V$  if and only if
the corresponding function $v=\left\langle V,x \right\rangle$ on the boundary of  $\Omega$ coincides with first order spherical harmonic.
\end{lmm}
This lemma is proved in Section \ref{harmonic}.
 
Based on the Lemma \ref{lem1} the proof of Theorem \ref{stbl} is done in Section \ref{thrm} and consists of two steps.
First, we consider special choice of vector fields $V$ when all domains $\Omega_t^V$ has the same
center of mass $x=0$ (Lemma \ref{lem2}).
Second, the case when center of mass moves, is
reduced to the previous step by suitable choice of non--autonomous vector field.

The proof of the Proposition \ref{pompeiu} about the Pompeiu problem equivalence to the shape optimization problem is 
done in Section \ref{secpmp}.

In Section \ref{sec5} the remark on generalization of Theorem \ref{stbl} to other domain functionals 
is done. Then a few numerical examples of shape evolution to the solution of  the Pompeiu  problem are presented.

\section{Shape derivatives of the domain functional }\label{shpdrv}

\begin{figure}
\begin{center}
\includegraphics[width=0.30\textwidth]{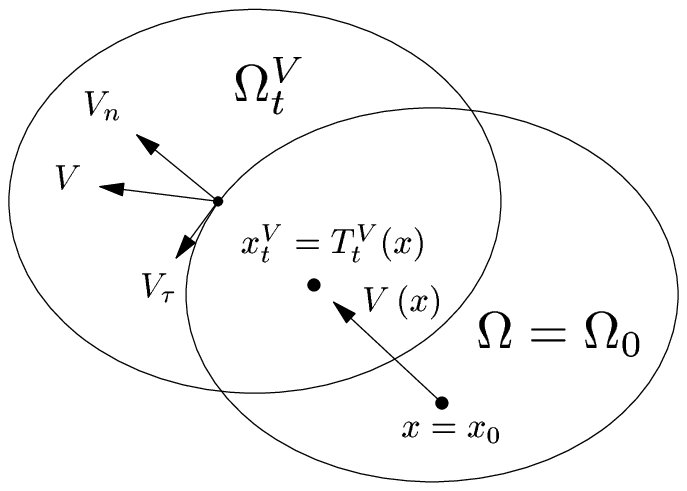}
\end{center}
\caption{}
\label{fig1}
\end{figure}

Let  $f:\left( 0, \infty \right)\rightarrow \xR$ be the smooth function. We
consider the domain functional
\begin{equation}
\it{F}\left[ \Omega \right] =\int\limits_\Omega \int\limits_\Omega
f\left( \left| x-y\right| \right) \xdif x\xdif y.
\end{equation}
Then, according to notions on page~\pageref{vect} (see also Fig.~\ref{fig1}),

\begin{equation}
\it{F}\left[ \Omega _t^V\right] =\int\limits_{\Omega
_0}\int\limits_{\Omega _0}f\left( \left| T_t^V\left( x\right) -T_t^V\left(
y\right) \right| \right) \det DT_t^V\left( x\right) \det DT_t^V\left(
y\right) \xdif x\xdif y,
\end{equation}
here $\det DT_t^V$ denotes the Jacobian of the mapping $T_t^V:\xR \times \xR^n\rightarrow \xR^n$,
$t$ considered as parameter.

Using equalities
\begin{gather}
\frac d{dt}T_t^V\left( x\right) =V\left( T_t^V\left( x\right) \right) ,\\
T_0^V\left( x\right) =x,\\
\frac d{dt}|_{t=0}\det DT_t^V\left( x\right) =\div V\left( x\right) ,
\end{gather}
(see \cite{Zol} on calculus technique) and verifying the equality
\begin{multline}
\frac d{dt}|_{t=0}f\left( \left| T_t^V\left( x\right) -T_t^V\left( y\right)
\right| \right) =\\
\left\langle\nabla _xf\left( \left| x-y\right| \right) ,V\left( x\right)
\right\rangle +\left\langle \nabla _yf\left( \left| x-y\right| \right) ,V\left( y\right)
\right\rangle 
\end{multline}
one gets
\small
\begin{equation}\label{eq5} 
\begin{split}
\frac d{dt}|_{t=0}\it{F}\left[ \Omega _t^V\right] =
\int\limits_{\Omega _0}\int\limits_{\Omega _0}\{\left[ \left\langle \nabla
_xf\left( \left| x-y\right| \right) ,V\left( x\right) \right\rangle +\left\langle\nabla
_yf\left( \left| x-y\right| \right) ,V\left( y\right) \right\rangle \right] \\
+f\left( \left| x-y\right| \right) \left[\div V\left( x\right) +\div V\left(
y\right) \right] \}\xdif x\xdif y=\\
2\int\limits_{\Omega _0}\int\limits_{\Omega _0}\div_x\left[ f\left( \left|
x-y\right| \right) V\left( x\right) \right] \xdif x\xdif y=
2\int\limits_{\partial \Omega _0}\left[ \int\limits_{\Omega _0}f\left(
\left| x-y\right| \right) \xdif y\right] \left\langle V,\eta \right\rangle \left( x\right)
\xdif s_x,
\end{split}
\end{equation}
\normalsize
here $\eta $ is exterior unit normal vector field on $\partial \Omega _0.$

If $\Omega _0$ gives an extremum to $\it{F}$, then necessarily
\begin{equation}\label{eq6}
\int\limits_{\Omega _0}f\left( \left| x-y\right| \right) \xdif y=0
\end{equation}
for all $x\in \partial \Omega _0.$

Replacement of the ''initial moment'' $t=0$ by the arbitrary moment $t$ leads
to equality
\begin{equation}\label{eq7}
\frac d{dt}\it{F}\left[ \Omega _t^V\right] =2\int\limits_{\partial
\Omega _t^V}\left[ \int\limits_{\Omega _t^V}f\left( \left| x-y\right|
\right) \xdif y\right] \left\langle V,\eta _t^V\right\rangle \left( x\right) \xdif s_x,
\end{equation}
here $\eta _t^V$ denotes exterior unit normal vector field on $\partial \Omega
_t^V.$

Now let us calculate the second derivative
\begin{equation}
\it{\ddot F}\left[ \Omega _0;V,V\right] :=\frac{d^2}{dt^2}|_{t=0}\it{F}\left[ \Omega _t^V\right]
\end{equation}
in case when $\Omega _0$ is critical domain to functional $\it{F}$.

Applying formula (5.9) in page 1105 of  \cite{Zol} and because of  (\ref{eq6}) one gets
\small{
\begin{equation}\label{eq8}
\begin{split}
\frac 12\it{\ddot F}\left[ \Omega _0;V,V\right] =\frac d{dt}|_{t=0}\left\{
\int\limits_{\partial \Omega _t^V}\left[ \int\limits_{\Omega _t^V}f\left(
\left| x-y\right| \right) \xdif y\right] \left\langle V,\eta _t^V\right\rangle \left(
x\right) \xdif s_x\right\} =\\
\int\limits_{\partial \Omega _0}\{\frac d{dt}|_{t=0}\left[
\int\limits_{\Omega _0}f\left( \left| T_t^V\left( x\right) -T_t^V(y)\right|
\right) \det DT_t^V(y)\xdif y\right] \left\langle V,\eta _t^V\right\rangle \left( T_t^V\left(
x\right) \right) \\
+\left( n-1\right) \it{K}\left( x\right) \int\limits_{\Omega _0}f\left(
\left| x-y\right| \right) \xdif y\cdot \left\langle V,\eta \right\rangle ^2\left( x\right)
\}\xdif s_x=\\
\int\limits_{\partial \Omega _0}\left\{ \int\limits_{\Omega _0}\left[
\div_y\left( f\left( \left| x-y\right| \right) V\left(y\right)\right) + 
\left\langle \nabla _xf\left( \left| x-y\right| \right) ,V\left(
x\right) \right\rangle \right]\xdif y \right\} \left\langle V,\eta \right\rangle \left( x\right) \xdif s_x=\\
\int\limits_{\partial \Omega _0}\int\limits_{\partial \Omega _0}f\left(
\left| x-y\right| \right) \left\langle V,\eta \right\rangle \left( x\right) \left\langle
V,\eta \right\rangle (y)\xdif s_x\xdif s_y\\
+\int\limits_{\partial \Omega _0}\left\langle \nabla _x\int\limits_{\Omega
_0}f\left( \left| x-y\right| \right) \xdif y,V\left( x\right) \right\rangle \left\langle V,\eta \right\rangle \left( x\right) \xdif s_x,
\end{split}
\end{equation}
}
\normalsize
here $\it{K}\left( x\right)$ denotes mean curvature of $\partial \Omega _0$ at the point $x$.

If function  $f\left( \left| \cdot \right| \right) $ is positive definite and $\it{F}\left[ \Omega_0 \right] =0$,
then indicator function $\chi _{\Omega _0}$ minimizes functional
\begin{equation}
\it{F}\left[ \varphi \right] =\int \int f\left( \left| x-y\right| \right)
\varphi \left( x\right) \varphi (y)\xdif x\xdif y
\end{equation}
defined on the set of compactly supported functions. Then necessarily
\begin{equation}
\int\limits_{\Omega _0}f\left( \left| x-y\right| \right) \xdif y=0
\end{equation}
for all $x\in \xR^n.$ As a consequence, the second term
in the last part of (\ref{eq8}) is zero.

So if $\it{F}\left[ \Omega_0 \right] =0$, then
\begin{equation}\label{eq12}
\it{\ddot F}\left[ \Omega _0;V,V\right] =2\int\limits_{\partial \Omega
_0}\int\limits_{\partial \Omega _0}f\left( \left| x-y\right| \right) v\left(
x\right) v(y)\xdif s_x\xdif s_y,
\end{equation}
for all $V\in C^1\left( \xR^n,\xR^n\right) .$ Here
\begin{equation}
v\left( x\right) =\left\langle V,\eta \right\rangle \left( x\right)
\end{equation}
denotes the normal component of vector field $V$ on $\partial \Omega _0.$

\section{Proof of Lemma \ref{lem1}}\label{harmonic}

The proof is based on standard properties of Fourier transform and
of Bessel functions.  For them we refer to \cite{StW} and \cite{Wat}.

Let $B$ be a ball of radius $R$ centered at $x=0$ and $\it{F}\left[ B\right]=0$ for the domain functional (\ref{funk}) with
positive definite function $f\left( \left| \cdot\right| \right)$. Then 
\begin{equation}\label{eq9}
\begin{split}
0=\int\limits_B\int\limits_Bf\left( \left| x-y\right| \right) \xdif x\xdif y=
\int\limits_0^\infty \int\limits_{\left| \xi \right| =r}\left| \hat \chi
_B\right| ^2\left( \xi \right) \xdif s_\xi \xdif \mu \left( r\right) =\\
\left( 2\pi
R\right) ^n\int\limits_0^\infty \int\limits_{\left| \xi \right| =r}\frac{J_{\frac n2}^2\left( R\left| \xi \right| \right) }{\left| \xi \right| ^n}\xdif s_\xi \xdif \mu \left( r\right) 
=\left( 2\pi R\right) ^n\omega _{n-1}\int\limits_0^\infty r^{-1}J_{\frac
n2}^2\left( Rr\right) \xdif \mu \left( r\right),
\end{split}
\end{equation}
here $\omega _{n-1}$ denotes the surface measure of the unit sphere.

Because function $\frac{J_{\frac n2}\left(t\right) }{t^{\frac{^n}2}}$ is entire and $\mu$ is positive Borel measure 
of polynomial growth \cite{RS}, equality in (\ref{eq9}) is possible 
only if measure $\mu$ is discrete and
\begin{equation}\label{eq10}
\supp\mu \subset \left\{ t|J_{\frac n2}\left( Rt\right) =0\right\} .
\end{equation}
Now, let us consider equality $\it{\ddot F}\left[ B;V,V\right] = 0$. We have
\begin{equation}
\begin{split}
\it{\ddot F}\left[ B;V,V\right] =
\int\limits_{\partial B}\int\limits_{\partial B}f\left( \left| x-y\right|
\right) v\left( x\right) v\left( y\right) \xdif s_x\xdif s_y=\\
\int\limits_0^\infty \int\limits_{\left| \xi \right| =r}\left| \left( \chi
_{\partial B}v\right) ^{\wedge }\right| ^2\left( \xi \right) \xdif s_\xi \xdif \mu
\left( r\right).
\end{split}
\end{equation}
Function $v\in L_2\left( \partial B\right) $ enables decomposition
\begin{equation}
v\left( R\xi \right) =\sum\limits_{k=0}^\infty c_kY_k\left( \xi \right)
,\quad\left| \xi \right| =1,
\end{equation}
for spherical harmonic $Y_k$ of order $k$ \cite{StW}.

Then
\begin{equation}
\begin{split}
\it{\ddot F}\left[ B;V,V\right] =
\int\limits_0^\infty r^{n-1}\int\limits_{\left| \xi \right|
=1}\sum\limits_{k=0}^\infty \left| c_k\right| ^2\left| \left( \chi
_{\partial B}Y_k\right) ^{\wedge }\right| ^2\left( \xi r\right) \xdif s_\xi \xdif \mu
\left( r\right) =\\
=\int\limits_0^\infty r^{n-1}\int\limits_{\left| \xi \right|
=1}\sum\limits_{k=0}^\infty \left| a_k\right| ^2\frac{J_{\frac{n+2k-2}2}^2\left( Rr\right) }{r^{n-2}}Y_k^2\left( Rr\right) \xdif s_\xi \xdif \mu \left(
r\right),
\end{split}
\end{equation}
since
\begin{equation}
\int\limits_{\left| \xi \right| =1}Y_k\left( \xi \right) e^{irR\left\langle \xi
,x\right\rangle }\xdif s_\xi =const\cdot \frac{J_{\frac{n+2k-2}2}\left( Rr\right) }{r^{\frac{^{n-2}}2}}Y_k\left( x\right)
\end{equation}
and $const\neq 0$ \cite{StW}.

Bessel functions of different order has nonintersecting sets of zeroes
\cite{Wat}, so due to relation (\ref{eq10}) equality
\begin{equation}
\int\limits_{\partial B}\int\limits_{\partial B}f\left( \left| x-y\right|
\right) v\left( x\right) v\left( y\right) \xdif s_x\xdif s_y=0
\end{equation}
is valid if and only if
\begin{equation}
\begin{split}
\int\limits_{\partial B}\int\limits_{\partial B}f\left( \left| x-y\right|
\right) v\left( x\right) v\left( y\right) \xdif s_x\xdif s_y=
\left| a_1\right|
^2\omega _{n-1}\int\limits_0^\infty rJ_{\frac n2}^2\left( Rr\right) \xdif \mu
\left( r\right) =\\
\int\limits_{\partial B}\int\limits_{\partial B}f\left( \left| x-y\right|
\right) Y_1\left( x\right) Y_1\left( y\right) \xdif s_x\xdif s_y,
\end{split}
\end{equation}

where $Y_1$ is the first order spherical harmonic.

\begin{rmrk}\label{rmrk2}
Function $Y_1$ has the form
\begin{equation}
Y_1\left( x\right) =b_1x_1+b_2x_2+...+b_nx_n=\left\langle b,x\right\rangle ,
\end{equation}
\end{rmrk}
i.e. on the surface of the ball $B$ it coincides with a normal component of a
constant vector field.
\clearpage

\section{Stability of the ball}\label{thrm}

In this Section we shall prove Theorem \ref{stbl}.

Because of Lemma \ref{lem1}, Remarks \ref{rmrk1} and \ref{rmrk2} it is sufficient to verify stability of 
the ball $\Omega =\Omega _0$ centered at $x=0$ in direction of vector field $V$ when function 
$v\left( x\right) =\left\langle V\left( x\right) ,x\right\rangle $ on the boundary $\partial \Omega$ is the
first order spherical harmonic.
\begin{lmm}\label{lem2}
Suppose that point $x=0$ is mass center of the domain $\Omega_t^V$ for all $t$ close enough to zero and
function $v\left( x\right) =\left\langle V\left( x\right) ,x\right\rangle $ on the boundary of the ball $\Omega_0$
is the first order spherical harmonic . Then $v=0$ on $\partial \Omega_0$.
\end{lmm}
\begin{proof}
For all $t$ close enough to zero and for any $i=\overline{1,n}$ equalities
$$
\int\limits_{\Omega _t^V}x_i\xdif x=0
$$
and
$$
0=\frac d{dt}|_{t=0}\int\limits_{\Omega _t^V}x_i\xdif x=\int\limits_{\partial \Omega _0}x_iv\left( x\right) \xdif s
$$
holds. Consequently, $v=0$ on the boundary $\partial \Omega_0$. 
\end{proof}
In general case, mass center $\bar x_t$ of the domain $\Omega _t^V$ is not zero. We consider the flow 
$\tilde T_t^V\left( x\right) :=T_t^V\left( x\right) -\bar x_t$. This flow
generates the family of domains $\tilde \Omega _t^V$ which coincide as rigid bodies with $\Omega _t^V$ 
and are centered at $x=0$. 
Of course, the new vector field 
$$
\tilde V\left( t,x\right) :=\frac d{dt} \tilde T_t^V\left( x\right)
$$
is non--autonomous. Fortunately, shape derivative method is applied also in the case of  non--autonomous vector 
fields $V\in \xCzero\left( I,\xCone\left( \xR^n,\xR^n\right) \right) $,
$I\subset \xR$, $0\in I$ \cite{Zol}. Moreover, withdrawal of formula (\ref{eq12}) in Section \ref{shpdrv} remains the same in the case of 
$v\left( x\right) =\left\langle \tilde V\left( 0,x\right), \eta \left( x\right) \right\rangle $.

\section{Shape optimization and the Pompeiu problem}\label{secpmp}

In this Section we shall prove Proposition \ref{pompeiu}.

{\bf Proof of a)}\quad Let equation (\ref {eq1}) has bounded domain $\Omega$ as its own solution for
some $f\neq 0$. By Bochner-Schwartz theorem \cite{RS} and
because of spherical symmetry of f$\left( \left| \cdot \right| \right) ,$
there exists $\mu $  --- positive Borel measure on $\left( 0, \infty \right)$ of polynomial
growth:
\begin{equation}
\begin{split}
\it{F}\left[ \Omega\right] =\int\limits_{\Omega}\int\limits_{\Omega}f\left( \left| x-y\right| \right) \xdif x\xdif y=
\int\limits_{\xR^n}\left| \hat \chi _{\Omega}\right| ^2\left( \xi
\right) \xdif \mu \left( \left| \xi \right| \right) =\\ \int\limits_0^\infty
\int\limits_{\left| \xi \right| =r}\left| \hat \chi _{\Omega}\right|
^2\left( \xi \right) \xdif s_\xi \xdif \mu \left( r\right) ,
\end{split}
\end{equation}
where
\begin{equation}
\hat \chi _{\Omega}\left( \xi \right) =\int\limits_\Omega e^{i\left\langle \xi ,x\right\rangle }\xdif x.
\end{equation}
Then equality
\begin{equation}
\it{F}\left[ \Omega\right] =0
\end{equation}
implies the existence of $\lambda \geq 0:$
\begin{equation}\label{eq3}
\int\limits_{\left| \xi \right| =\lambda }\left| \hat \chi _{\Omega}\right| ^2\left( \xi \right) \xdif s_\xi =0.
\end{equation}
Consequently,
\begin{equation}\label{eq4}
\hat \chi _{\Omega}\left( \xi \right) =0
\end{equation}
for all $\left| \xi \right| =\lambda .$
Because of $\hat \chi _{\Omega}\left( 0\right) =\mes\Omega>0,$ one has $\lambda >0.$

{\bf Proof of b)}\quad Suppose that $\Omega$ has no the Pompeiu property. Then (\cite{Ber1},
\cite{BST}) there exists $\lambda >0$ equality (\ref {eq4}) is true and, therefore,
\begin{equation}
\begin{split}
0=\int\limits_{\left| \xi \right| =\lambda }\left| \hat \chi _{\Omega}\right| ^2\left( \xi \right) \xdif s_\xi =
\int\limits_{\left| \xi \right| =
\lambda }\left[ \int\limits_{\Omega}\int\limits_{\Omega}e^{i\left\langle \xi ,x-y\right\rangle }\xdif x\xdif y\right] \xdif s_\xi =\\
\int\limits_{\Omega}\int\limits_{\Omega}\left[ \int\limits_
{\left|\xi \right| =\lambda }e^{i\left\langle \xi ,x-y\right\rangle }\xdif s_\xi \right] \xdif x\xdif y=
\left( 2\pi \lambda \right) ^{\frac n2}\int\limits_{\Omega}\int\limits_{\Omega}\frac{J_{\frac{n-2}2}\left( \lambda \left|
x-y\right| \right) }{\left| x-y\right| ^{\frac{n-2}2}}\xdif x\xdif y,
\end{split}
\end{equation}
here $J_p$ denotes the Bessel function.

This will end the proof of Proposition \ref{pompeiu}.

\section{Concluding remarks}\label{sec5}
{\bf 1.} Let $\it{F}$ be a domain functional defined on the set of bounded domains $\Omega \subset \xR^n$
and invariant under the rigid motions of $\xR^n$ .
Suppose that for the ball $B$ centered at $x=0$ the condition
\begin{itemize}
\item[(*)] \textit{if} \quad $\it{\ddot F}\left [ B;V,V\right ]=0$, 
\textit{then function} \quad
$v\left( x\right)=\left\langle V\left( x\right), x\right\rangle$ \textit{is first order spherical harmonic}
\end{itemize}
holds.
\begin{prpstn}\label{gen}
If $\it{\ddot F}\left [ B;V,V\right ]\geq 0$ for all $V\in \xCone\left( \xR^n,\xR^n\right)$ and the condition {\rm (*)}
holds, then the ball $B$ is weakly stable minimum point of domain functional $\it{F}$ in sense of Definition \ref{defstb}.
\end{prpstn}
\bigskip

{\bf 2.} Concerning the Schiffer problem. The result of Theorem \ref{stbl} does not cover ones stated in
\cite{AgS, Kob}. This is because spectral parameter $\lambda$ in (\ref{schif}) may vary with the domain $\Omega$.
Though this more general situation may be considered on the basis of
Lemma \ref{lem1} the proof requires additional techniques.
\bigskip

\begin{figure}
\includegraphics[width=0.30\textwidth]{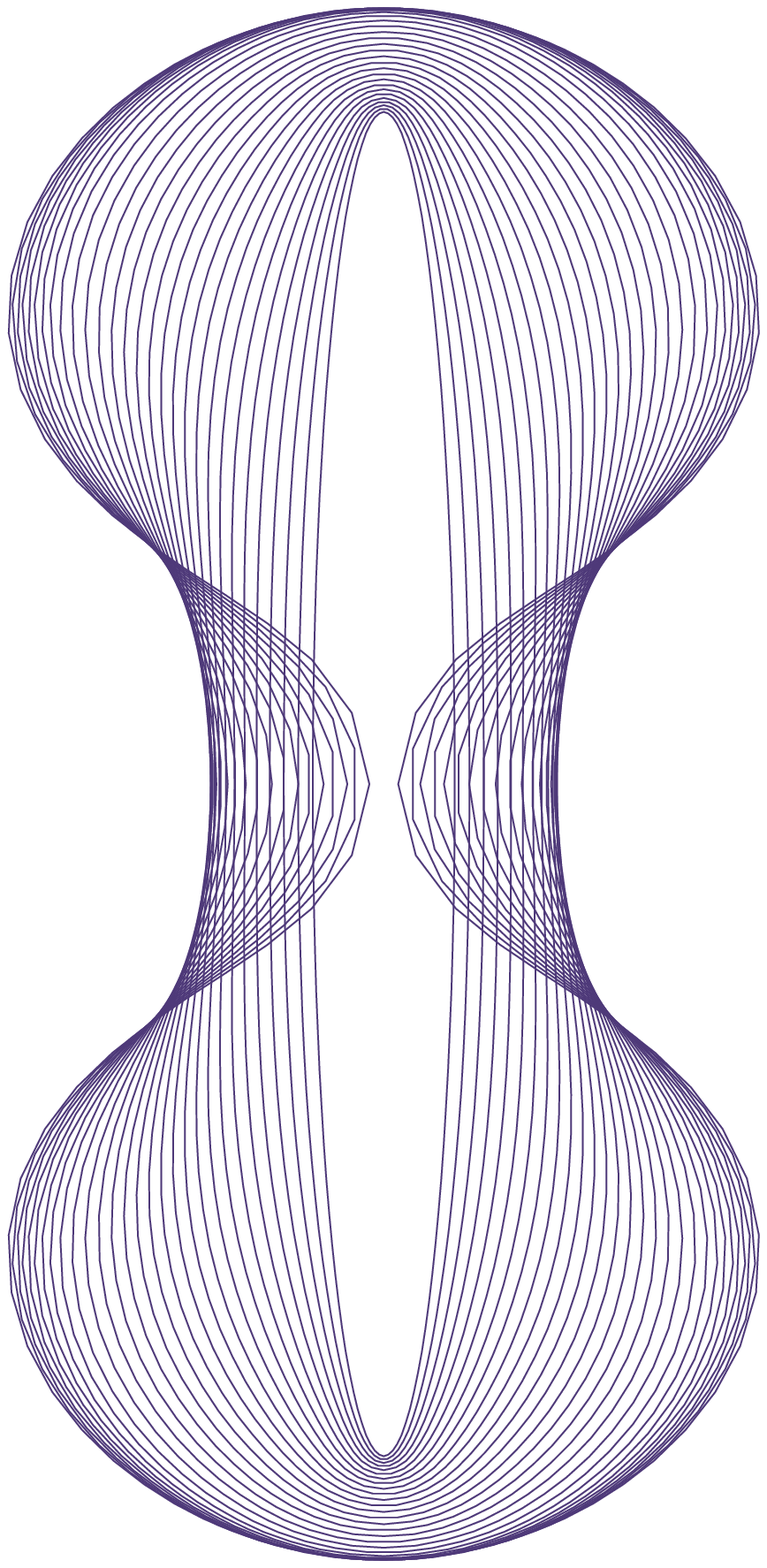}
\qquad
\includegraphics[width=0.30\textwidth]{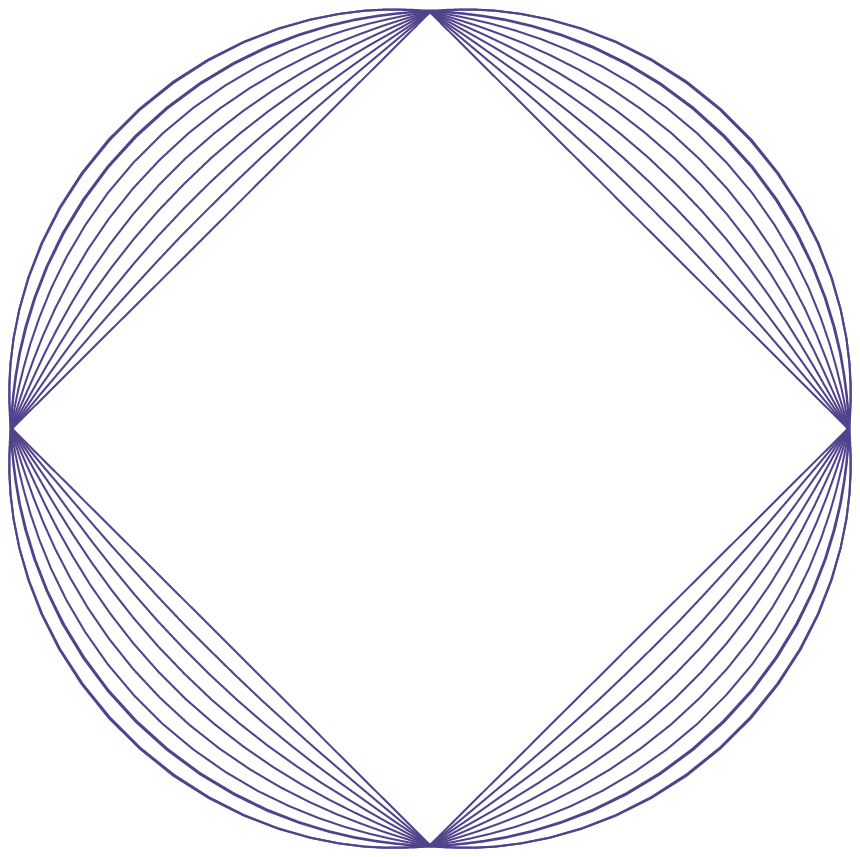}
\caption{}
\label{fig2}
\end{figure}

{\bf 3.} Methods based on the shape derivative allow numerical experiments in the Pompeiu problem.
In Fig.~\ref{fig2} evolution of the long thin ellipse and of the quadrate to the "figure" without the Pompeiu
property in the shape antigradient direction is presented.

\end{document}